\documentclass{article}

\usepackage{hyperref} 
\usepackage{amsfonts}

\newcommand{\dd}{{\, {\mathrm d}}}

\newtheorem{remark}{Remark}[section]
\newtheorem{theorem}{Theorem}[section]

\newtheorem{coroll}[theorem]{Corollary}

\begin{document}

\title{THE MEAN VALUE THEOREMS
AND A NAGUMO-TYPE UNIQUENESS THEOREM
FOR CAPUTO'S FRACTIONAL CALCULUS\\
(Corrected Version)}
\author{Kai Diethelm\footnote{Institut Computational Mathematics,
    Technische Universit\"at Braunschweig,
    Universit\"atsplatz 2,
    38106 Braunschweig, Germany,
    {\tt k.diethelm@tu-bs.de}}\ \footnote{Gesellschaft f\"ur numerische Simulation mbH,
    Am Gau\ss berg 2,
    38114 Braunschweig, Germany,
    {\tt diethelm@gns-mbh.com}}
}

\maketitle

\begin{center}
\emph{Dedicated to the memory of my teacher, Professor 
Dr.\ Helmut Bra\ss.}
\end{center}

\begin{abstract}
We generalize the classical mean value theorem of differential 
calculus by allowing the use of a Caputo-type fractional 
derivative instead of the commonly used first-order derivative.
Similarly, we generalize the classical mean value theorem 
for integrals by allowing the corresponding fractional integral,
viz.\ the Riemann-Liouville operator,
instead of a classical (first-order) integral.
As an application of the former result we then prove a uniqueness
theorem for initial value problems involving Caputo-type fractional
differential operators. This theorem generalizes the classical 
Nagumo theorem for first-order differential equations.

\smallskip

\framebox{\parbox{0.8\textwidth}{The original version of this paper, published in 
\emph{Fract.\ Calc.\ Appl.\ Anal.} 15 (2012), pp.\ 304--313,
unfortunately contained an error in Corollary \ref{cor:mvti-f}
that was then carried forward to the later parts of the paper.
This version contains the corrected form of the document.}}

\medskip

{\it MSC 2010\/}: Primary 26A33; Secondary 34A08, 34A12

\medskip

{\it Key Words and Phrases}: Caputo derivative, 
mean value theorem of differential calculus,
mean value theorem of integral calculus,
fractional differential equation,
Nagumo's theorem,
uniqueness of solutions

\end{abstract}

\section{Introduction}

In this short note we shall demonstrate that two well known results connected
to classical analysis, namely the mean value theorems of differential and of
integral calculus, can be extended to fractional 
calculus, i.e.\ they can be generalized by replacing the first derivatives and
integrals, respectively, by derivatives and integrals of non-integer order. 
As an application of our results we shall then prove a new uniqueness theorem
for a class of fractional differential equations.

It is well known \cite{Di,KST,SKM} that  
many different concepts for fractional derivatives exist. The appropriate
approaches for our purposes are the Riemann-Liouville integral and the 
Caputo derivative of order $\alpha > 0$
with starting point $a \in \mathbb R$ \cite[Chapters 2 and 3]{Di}.
The latter is denoted and defined by  
\begin{equation}
\label{eq:cap-der}
D_{*a}^\alpha f := D_a^\alpha[ f - T_{\lceil \alpha \rceil-1}[f;a]]
\end{equation}
where $\lceil \cdot \rceil$ denotes the ceiling function that rounds up to the
nearest integer, $T_m[f;a]$ is the $m$th degree Taylor polynomial for the
function $f$, centered at $a$, and $D_a^\alpha$ is the fractional
differential operator of Riemann-Liouville type, defined by
\begin{equation}
\label{eq:rl-der}
D_a^\alpha f := D^{\lceil \alpha \rceil} J_a^{\lceil \alpha \rceil - \alpha} f,
\end{equation}
with $D^m$ being the classical differential operator of integer order $m$ and
$J_a^\beta$ denoting the Riemann-Liouville integral operator of order $\beta
\ge 0$, given by
\begin{equation}
\label{eq:rl-int}
J_a^\beta f(x) := 
\cases{
         f(x) & for $ \beta = 0$, \cr
     \displaystyle \frac1{\Gamma(\beta)} \int_a^x (x-t)^{\beta-1} f(t) \dd t
              & for $\beta > 0$. \cr
}
\end{equation}
for $x \ge a$. Under appropriate
smoothness assumptions on $f$, we may write
\begin{equation}
\label{eq:cap-der-2}
D_{*a}^\alpha f :=  J_a^{\lceil \alpha \rceil - \alpha} D^{\lceil \alpha \rceil} f.
\end{equation}
A sufficient condition for this identity to hold is, e.g., that $f$ has an
absolutely continuous $\lceil \alpha - 1\rceil$st derivative \cite[Theorem 3.1]{Di}.

We note that the Riemann-Liouville derivatives defined in eq.\ (\ref{eq:rl-der}) 
have been investigated and used in mathematical analysis for a very long time; 
cf.\ the exhaustive treatment in the classical monograph \cite{SKM}. However, 
in the last couple of decades, the differential operators of Caputo's type as
given in eq.\ (\ref{eq:cap-der}) have been found to be very important because 
they have proven to be highly useful for the mathematical modeling of many
phenomena of great interest in areas like mechanics, life sciences, finance,
etc; cf.\ the recent monograph \cite{Di} and the references cited therein.
Therefore it is an important task to complete the investigation of the analytical 
properties of the Caputo operators.

It is very natural to combine the Caputo derivative with the
Riemann-Liouville integral because these operators are one-sided inverses of
each other. This is clearly demonstrated by the fractional version of the 
fundamental theorem of calculus, 
\begin{equation}
\label{eq:fund}
D_{*a}^\alpha J^\alpha_a f(x) = f(x)
\end{equation}
for $f \in C[a,b]$ and $\alpha > 0$ (see, e.g., \cite[Theorem 3.7]{Di}).

For further information about the analytical properties of all these
operators, their relations to other fractional derivatives, and their many
applications in various mathematical, technical and other scientific
disciplines, we refer to the monographs of Samko et al.\ \cite{SKM}, Kilbas et
al.\ \cite{KST}, and Diethelm \cite{Di}.

\section{The fractional versions of the mean value theorems}
\label{sec:mvt}
\setcounter{equation}{0}

We shall now turn our attention towards our intended generalization of the two
above mentioned classical mean value theorems.

\subsection{Fractional mean value theorems of integral calculus}

First we look at the generalized mean value theorem of integral calculus. It reads
\begin{equation}
\label{eq:mvti}
\int_a^b f(t) g(t) \dd t = f(\xi) \int_a^b g(t) dt
\end{equation}
with some $\xi \in (a,b)$ if $f \in C[a,b]$, $g$ is Lebesgue integrable on
$[a,b]$ and $g$ does not have a change of sign in $[a,b]$ (cf., e.g., 
\cite[Theorem 85.6]{Heuser}). 

\begin{remark}
As a matter of fact, this result is stated in \cite{Heuser} and in some other 
textbooks under the assumption that $f$ is integrable in Riemann's rather than 
Lebesgue's sense; however, the method of proof can be used under the Lebesgue 
integrability assumption --- which is more appropriate for the fractional 
generalization that we have in mind --- in the same way (employing, in particular,
\cite[Theorem 36.4]{Heuser}).
\end{remark}
\begin{remark}
Many textbooks (e.g., \cite[\S 121, eq.\ (2)]{Brand}) provide
an even weaker formulation in the sense that they only prove that $\xi \in 
[a,b]$. However, for our applications in Section \ref{sec:nagumo} below, it
is important to know that a suitable $\xi$ can be found in the \emph{open} 
interval $(a,b)$. Therefore we have introduced Heuser's form of the result here.
\end{remark}

The fractional version of this theorem can be formulated as follows.

\begin{theorem}
\label{thm:mvti-f}
Let $a<b$, $\alpha > 0$ and $f \in C[a,b]$. Moreover, assume that $g$ is
Lebesgue integrable on $[a,b]$ and that $g$ does not change its sign in this
interval. Then, for almost every $x \in (a,b]$ there exists some $\xi \in (a,x) 
\subset (a,b)$ such that 
$$ J_a^\alpha (fg)(x) = f(\xi) J_a^\alpha g(x). $$
If additionally  
$\alpha \ge 1$ or $g \in C[a,b]$, then this result holds for every $x \in (a,b]$.
\end{theorem}

Clearly, the classical form (\ref{eq:mvti}) follows by setting $\alpha = 1$ and $x=b$.
\medskip

\emph{Proof.}
Under the assumptions of the theorem, we find that
$$
J_a^\alpha (fg)(x) = \frac1{\Gamma(\alpha)} \int_a^x (x-t)^{\alpha-1} f(t) g(t) \dd t.
$$
Now let us first assume that $\alpha \ge 1$ and $x \in (a,b]$. 
Then, the first factor of the integrand on the right-hand side of this
equation is continuous. Therefore, 
the function $\tilde g$ defined by $\tilde g(t) = (x-t)^{\alpha-1} g(t) / \Gamma(\alpha)$ is
integrable on $[a,x]$. Moreover, it does not change its sign in this
interval. Thus, by the classical result (\ref{eq:mvti}),
$$
J_a^\alpha (fg)(x) 
= \int_a^x f(t) \tilde g(t) \dd t 
= f(\xi) \int_a^x \tilde g(t) \dd t
= f(\xi) J_a^\alpha g(x).
$$
In the case that $0 < \alpha < 1$ and $g$ is continuous, the same line of proof 
works. Finally, if $0 < \alpha < 1$ and $g$ is only integrable, then we 
can still argue in a similar way, but the
integrability of $\tilde g$ holds only for almost all $x$ 
\cite[Theorem 4.2(d)]{Will}.
$\Box$\medskip

If we set $g(x) = 1$ in Theorem \ref{thm:mvti-f} and observe that
$J_a^\alpha g(x) = (x-a)^\alpha / \Gamma(\alpha+1)$, we immediately obtain the 
following corollary for $x=b$.

\begin{coroll}
\label{cor:mvti-f}
Let $a<b$, $\alpha > 0$ and $f \in C[a,b]$. Then, there exists some $\xi \in (a,b)$
such that
$$ J_a^\alpha f(b) = \frac 1 {\Gamma(\alpha+1)} (b-a)^\alpha f(\xi). $$
\end{coroll}

For $\alpha=1$ this reduces to the classical mean value theorem of integral calculus,
\begin{equation}
\label{eq:mvti2}
\int_a^b f(t) \dd t = (b-a) f(\xi),
\end{equation}
that, of course, can also be obtained from (\ref{eq:mvti}) by setting $g(x) = 1$ 
for all $x$.

\subsection{The fractional mean value theorem of differential calculus}

Now we turn our attention towards the mean value theorem of differential calculus which 
states that
\begin{equation}
\label{eq:mvt-class}
\frac{f(b) - f(a)}{b-a} = f'(\xi)
\end{equation}
with some $\xi \in (a,b)$ if $a<b$ and $f \in C^1(a,b) \cap C[a,b]$
(see, e.g., \cite[\S 57]{Brand}). Our goal is to prove the following
generalization of this well known result:

\begin{theorem}
\label{thm:mvtd-f}
Let $\alpha > 0$ and $a < b$, and assume $f \in C^{\lceil \alpha \rceil-1}[a,b]$ 
to be such that $D^\alpha_{*a}f \in C[a,b]$. Then,
there exists some $\xi$ in $(a,b)$ such that 
$$
\frac{f(b) - T_{\lceil \alpha \rceil-1}[f;a](b)}{(b-a)^\alpha}
= \frac1{\Gamma(\alpha+1)} D_{*a}^\alpha f(\xi).
$$
\end{theorem}

Evidently, if $0 < \alpha \le 1$ then the Taylor polynomial in Theorem 
\ref{thm:mvtd-f} reduces to $T_{\lceil \alpha \rceil-1}[f;a](b) 
= T_0[f;a](b) = f(a)$,
and hence the mean value theorem takes the following form.

\begin{coroll}
\label{cor:mvtd-f}
Let $0 < \alpha \le 1$, $a < b$ and $f \in C[a,b]$ be such that 
$D^\alpha_{*a}f \in C[a,b]$. Then, there exists some $\xi \in (a,b)$ such that
$$
\frac{f(b) - f(a)}{(b-a)^\alpha}
= \frac 1{\Gamma(\alpha+1)} D_{*a}^\alpha f(\xi).
$$
\end{coroll}

This special case of our general result has already been derived in 
\cite[Theorem 1]{OS}.
Clearly, we recover the classical result (\ref{eq:mvt-class}) in the case 
$\alpha = 1$.

For the proof of Theorem \ref{thm:mvtd-f} we need a fractional version of the
Taylor expansion with integral representation of the remainder term whose
structure differs significantly from the result of \cite[Theorem 3]{OS}:

\begin{theorem}
\label{thm:taylor}
Let $\alpha > 0$ and $a < b$, and assume $f \in C^{\lceil \alpha \rceil-1}[a,b]$ 
to be such that $D^\alpha_{*a}f \in C[a,b]$. Then, for all $x \in [a,b]$,
\begin{equation}
\label{eq:taylor}
f(x) - T_{\lceil \alpha \rceil-1}[f;a](x)
= J_a^\alpha D_{*a}^\alpha f(x).
\end{equation}
\end{theorem}

Relations like (\ref{eq:taylor}) are well known in fractional calculus, 
cf.\ \cite[Corollary 3.9]{Di}. However, throughout the literature
these results are commonly proved under somewhat stronger 
assumptions on the function $f$ that would be too restrictive for the 
applications we have in mind. Therefore we now give a proof of 
eq.\ (\ref{eq:taylor}) under the weaker conditions mentioned above.

\medskip

\emph{Proof.}
Consider the functions $y_1:[a,b] \to \mathbb R$ and $y_2:[a,b] \to \mathbb R$
defined by $y_1(x) =  f(x) - T_{\lceil \alpha \rceil-1}[f;a](x)$ and $y_2(x)
:= J_a^\alpha D_{*a}^\alpha f(x)$. We need to prove that $y_1 = y_2$. To this
end we first note that $y_1$ and $y_2$ are continuous on $[a,b]$ by definition.
Moreover, 
\begin{equation}
\label{eq:taylor-ic}
\frac{\mathrm d^k}{\mathrm dx^k} y_j(0) = 0, \qquad k = 0,1,\ldots,\lceil
\alpha \rceil-1,  \quad j \in \{1,2\}.
\end{equation}
For $j=1$ this immediately follows from the well known properties of the
classical Taylor polynomial $T_{\lceil \alpha \rceil-1}[f;a]$; for $j=2$ it
follows from standard properties of the integral operator $J_a^\alpha$ in view
of the continuity of $D_{*a}^\alpha f$.

Next we note that 
\begin{equation}
\label{eq:taylor-de}
D^\alpha_{*a} y_j(x) = D_{*a}^\alpha f(x), \qquad j \in \{1,2\}.
\end{equation}
This follows for $j=1$ from eq.\ (\ref{eq:cap-der-2}) which implies that the
differential operator $D^\alpha_{*a}$ annihilates the Taylor polynomial, and
for $j=2$ it is a direct consequence of eq.\ (\ref{eq:fund}).

Thus, we have found that both functions $y_1$ and $y_2$ solve the
fractional-order initial value problem that consists of the differential
equation (\ref{eq:taylor-de}) and the initial conditions (\ref{eq:taylor-ic}),
and as the right-hand side of eq.\ (\ref{eq:taylor-de}) is continuous and
satisfies a Lipschitz condition with respect to $y_j$, we know from 
\cite[Theorem 2.2]{DF} (see also \cite[Theorem 6.5]{Di}) that this initial
value problem has a unique continuous solution. Therefore we have $y_1 = y_2$ as
desired.
$\Box$\medskip

\emph{Proof of Theorem \protect{\ref{thm:mvtd-f}}}
The fractional Taylor expansion
$$ 
f(b) - T_{\lceil \alpha \rceil-1}[f;a](b)
= J_a^\alpha D_{*a}^\alpha f(b)
$$
is valid because of Theorem \ref{thm:taylor}.
Moreover, the function $D_{*a}^\alpha f$ is continuous by assumption. 
We may thus apply Corollary \ref{cor:mvti-f}  to the right-hand side 
of this equation (with $D_{*a}^\alpha f$ taking the role of $f$),
and the claim follows immediately.
$\Box$\medskip

\section{A Nagumo-type uniqueness theorem for fractional-order initial value
  problems}
\label{sec:nagumo}
\setcounter{equation}{0}

The results proved in Section \ref{sec:mvt} enable us to prove a
generalization of Na\-gu\-mo's classical uniqueness theorem for first order
initial value problems that states that the initial value problem
\begin{equation}
\label{eq:ivp-1}
y'(x) = f(x, y(x)), \qquad y(0) = y_0, 
\end{equation}
has at most one solution on the interval $[0,b]$ if the function $f:[0,b] \times
\mathbb R \to \mathbb R$ on the right-hand side of the
differential equation is continuous at the initial point $(0, y_0)$ and 
satisfies the inequality
\begin{equation}
\label{eq:nagumo-condition1}
x | f(x, y_1) - f(x, y_2) | \le |y_1 - y_2|
\end{equation}
for all $x \in [0,b]$ and all $y_1, y_2 \in \mathbb R$ (see, e.g., Nagumo's
original work \cite{Nagumo} or the very elegant proof given by Diaz and Walter
\cite{DW} whose
path we shall follow below in the proof of Theorem \ref{thm:nagumo-f}). 

A number of generalizations of this result to fractional differential
equations containing Riemann-Liouville differential operators have been 
developed recently \cite{Baleanu1,LL}.
Our fractional version of this result, however, deals with differential
equations containing Caputo-type derivatives, for which no corresponding
result seems to be known so far. It complements the existing literature on uniqueness
results for such equations \cite[Chapter 6]{Di} and reads as follows.

\begin{theorem}
\label{thm:nagumo-f}
Let $\alpha \in (0,1)$, $b>0$ and $y_0 \in \mathbb R$.
If the function $f:[0,b] \times \mathbb R \to \mathbb R$ 
is continuous at $(0, y_0)$ and satisfies the inequality
\begin{equation}
\label{eq:nagumo-condition-frac}
x^\alpha | f(x, y_1) - f(x, y_2) | \le \Gamma(\alpha+1) |y_1 - y_2|
\end{equation}
for all $x \in [0,b]$ and all $y_1, y_2 \in \mathbb R$
then the initial value problem 
\begin{equation}
\label{eq:ivp-f}
D_{*0}^\alpha y(x) = f(x, y(x)), \qquad y(0) = y_0, 
\end{equation}
has at most one continuous solution $y$ on $[0,b]$ satisfying $D_{*0}^\alpha y
\in C[0,b]$.
\end{theorem}

It is immediately evident that the case $\alpha = 1$ reproduces the classical
result mentioned above. 

\medskip

\emph{Proof.}
  Assume that the initial value problem (\ref{eq:ivp-f}) has two continuous solutions $z$
  and $\tilde z$, say, on $[0,b]$. We then have to prove that $z = \tilde z$. To
  this end, we define the function $w$ on $[0,b]$ by 
  $$
   w(x) := \cases{
              x^{-\alpha} |z(x) - \tilde z(x)| & for $ x \in (0,b]$, \cr
                0                            & for $ x = 0$. \cr
           }
  $$
  It is then obvious that $w$ in continuous on
  $(0,b]$. Moreover, since both $z$ and $\tilde z$ solve the initial value
  problem and in view our Corollary \ref{cor:mvtd-f}, we find that
  \begin{eqnarray*}
    w(x) &=& x^{-\alpha} |z(x) - \tilde z(x)| 
          =  \left| x^{-\alpha} ([z(x) - z(0)] - [\tilde z(x) - \tilde z(0)]) \right| \\
         &=& \frac1{\Gamma(\alpha+1)} |D_{*0}^\alpha z(\xi) - D_{*0}^\alpha \tilde z(\xi)| 
          = \frac1{\Gamma(\alpha+1)} |f(\xi, z(\xi)) - f(\xi, \tilde z (\xi)) |
  \end{eqnarray*}
  with some $\xi \in (0,x)$. Thus, as $x \to 0$ we also have $\xi \to 0$, 
  and since $z$ and $\tilde z$ are continuous, 
  we also see that $z(\xi) \to z(0) = y_0$ and $\tilde z(\xi) \to \tilde z(0) = y_0$. 
  These observations, combined with the continuity
  assumption for $f$ at the point $(0, y_0)$, yield that $w(x) = f(\xi, z(\xi)) -
  f(\xi, \tilde z (\xi)) \to f(0, y_0) - f(0, y_0) = 0$ for $x
  \to 0$. In particular, $\lim_{x \to 0} w(x)$ exists and coincides with $w(0)$
  which demonstrates that $w$ is continuous also at the origin.

  Now we assume that $w \not \equiv 0$ on $[0,b]$ and define
  $$ \eta := \inf \left\{ x \in [0,b] : w(x) = \sup_{t \in [0,b]} w(t) \right\} . $$ 
  Because $w$ is continuous and nonnegative and $w(0) = 0$, we conclude that 
  $w(\eta) = \sup_{t \in [0,b]} w(t)$ and that 
  \begin{equation}
  \label{eq:proof-nagumo-f1}
  w(\tau) < w(\eta) \mbox{ for all } \tau \in [0,\eta).
  \end{equation}
  By another application of our mean value theorem of Caputo's fractional
  differential calculus (Corollary \ref{cor:mvtd-f}) and the Nagumo condition
  (\ref{eq:nagumo-condition-frac}), we derive
  \begin{eqnarray*}
    w(\eta) 
      &=& \eta^{-\alpha} |[z(\eta) - z(0)] - [\tilde z(\eta) - \tilde z(0)]| \\
      &=& \frac1{\Gamma(\alpha+1)} | D_{*0}^\alpha z(\tau) - D_{*0}^\alpha \tilde z(\tau)| \\
      &=& \frac1{\Gamma(\alpha+1)} |f(\tau, z(\tau)) - f(\tau, \tilde z(\tau))| 
      \le w(\tau)
  \end{eqnarray*}
  with some $\tau \in (0,\eta)$ which contradicts
  eq.\ (\ref{eq:proof-nagumo-f1}).
  It thus follows, as desired, that $w$ vanishes identically.
$\Box$\medskip

\begin{remark}
As in the integer order case \cite[p.\ 93]{DW}, we observe that the continuity
requirement for $f$ at the point $(0,y_0)$ that we had imposed in Theorem
\ref{thm:nagumo-f} is essential. Specifically, let $y_0 = 0$ and define
$$ 
f(x,y) = 
\cases{ 
                  \Gamma(\alpha+1) &  for  $y > x^\alpha$, \cr
      \Gamma(\alpha+1) x^{-\alpha} y &  for  $0 < y \le x^\alpha$, \cr
                                 0 &  for  $y \le 0$. \cr
}
$$
Condition (\ref{eq:nagumo-condition-frac}) is obviously satisfied. 
Clearly, the function $f$ is not continuous at $(0,y_0)$. Moreover, we see
that $y(x) = c x$ is a solution of the initial value problem (\ref{eq:ivp-f})
for all $c \in [0,1]$, so we have infinitely many solutions.
\end{remark}

A slight modification of the hypotheses leads us to an existence and
uniqueness theorem. 

\begin{theorem}
Let $\alpha \in (0,1)$, $b>0$ and $y_0 \in \mathbb R$.
If the function $f:[0,b] \times \mathbb R \to \mathbb R$ 
is continuous in its domain of definition and satisfies the inequality
(\ref{eq:nagumo-condition-frac})
for all $x \in [0,b]$ and all $y_1, y_2 \in \mathbb R$
then the initial value problem (\ref{eq:ivp-f})
has at exactly one continuous solution $y$ on $[0,b]$.
\end{theorem}

\emph{Proof.}
As we have now assumed the continuity of $f$ everywhere, the existence of a
solution is clear from the fractional version of Peano's existence theorem
(see \cite[Theorem 2.1]{DF} or \cite[Theorem 6.1]{Di}). Moreover,
the continuity of $f$ in combination with the differential equation of
(\ref{eq:ivp-f}) implies that each continuous solution $y$ also satisfies
$D^\alpha_{*0} y \in C[0,b]$. Therefore, the uniqueness of the continuous
solution follows from Theorem \ref{thm:nagumo-f}.
$\Box$\medskip

\begin{remark}
Recently, fractional terminal value problems (sometimes also called boundary
value problems), i.e.\ problems where the differential equation from 
eq.\ (\ref{eq:ivp-f}), viz.\ 
$$ D_{*0}^\alpha y(x) = f(x, y(x)) $$
for some $\alpha \in (0,1)$, is combined with a condition of the form
$$ y(b) = y_0 $$
for some $b > 0$, have drawn a lot of attention \cite{Di:tvp,DF:tvp,FM}. 
It would be of interest to find out whether the conditions of our Nagumo-type 
result, Theorem \ref{thm:nagumo-f}, can be modified to prove a corresponding
uniqueness theorem for this class of problems too. We intend to investigate
this question in a future paper. 
\end{remark}



\end{document}